\documentclass[review]{autart}

\usepackage{lineno,hyperref}
\modulolinenumbers[5]

\usepackage{amssymb}
\usepackage{amsmath}

\usepackage{graphicx}
\usepackage{psfrag}

\newcommand{\R}{\mathbb{R}}

\newcommand{\be}{\end{eqnarray*}}
\newcommand{\ee}{\end{eqnarray*}}
\newcommand{\ben}{\begin{eqnarray}}
\newcommand{\een}{\end{eqnarray}}
\newtheorem{lemma}{Lemma}[section]

\newtheorem{proposition}[lemma]{Proposition}

\newtheorem{definition}[lemma]{Definition}

\journal{Automatica}

\begin{document}
\begin{frontmatter}
\title{Decomposition of Differential Games}

\author[mymainaddress]{Adriano Festa}
\ead{festa@ensta.fr}
\author[mysecondaryaddress]{ Richard Vinter}
\ead{r.vinter@imperial.ac.uk}
\address[mymainaddress]{ENSTA ParisTech, 828, Boulevard des Maréchaux, 91120 Palaiseau, FR}
\address[mysecondaryaddress]{Imperial College, EEE Department, Exhibition Road, SW7 2BT  London, UK}

\begin{abstract}
This paper  provides a decomposition technique for the purpose of simplifying the solution of certain zero-sum differential games. The games considered terminate when the state reaches a target, which can be expressed as the union of a collection of target subsets; the decomposition consists of replacing the original target by each of the target subsets. The value of the original game is then obtained as the lower envelope of the  values of the collection of games resulting from the decomposition, which can be much easier to solve than the original game. Criteria are given for the validity of the decomposition. The paper includes examples, illustrating the application of the technique to pursuit/evasion games, where the decomposition arises from considering the interaction of individual pursuer/evader pairs.
\end{abstract}


\end{frontmatter}

\section{Introduction}

\label{section1}
We propose a decomposition technique to simplify the solution of  zero-sum differential games that involve two players (the $a$-player and the $b$-player), whose actions govern the evolution of the state $x$. The  state trajectory associated with open loop policies $a(.)$ and $b(.)$ (`open loop policies' are defined below), for a specified initial state $x_{0}$, is given by the (absolutely continuous) solution of the differential equation
$$
\left\{
\begin{array}{l}
\dot x(t)\; =\;f(x(t), a(t),b(t)),\; \mbox{ a.e.}
\\
x(0)=x_{0}\;.
\end{array}
\right.
$$
Here, $f(.,.,.): \R^{n} \times \R^{m_{1}}\times  \R^{m_{2}} \rightarrow \R^{n}$ is a given function. Open loop policies $a(.)$ and $b(.)$ of the two players  take values in specified sets $A \subset \R^{m_{1}}$ and $B \subset \R^{m_{2}}$  respectively. We write the solution $x(t; x_{0},a(.),b(.))$.  It is assumed that hypotheses are imposed on the data ensuring that a solution exists and it is unique. 
We also specify a closed set ${\mathcal{T}} \subset{\R^{n}}$ called the `target'.
The first entry time 
$\tau$ for $x(t; x_{0},a(.), b(.))$ is.
$$
\tau \,:=\, \sup \{t\,|\, x(t; x_{0}, a(.), b(.))  \notin {\mathcal{T}} \}\;.  
$$
%
%
%
%
Let ${\mathcal{A}}$ and ${\mathcal{B}}$ be the spaces of open loop policies for the $a$-player and $b$-player respectively, namely
\begin{eqnarray*}
{\mathcal{A}}&:=&\{ a(.):[0,\infty) \rightarrow \R^{m_{1}}\,|\, a(.) \mbox{ meas. and } a(t) \in A \mbox{ a.e. }  \},
\\
{\mathcal{B}}&:=&\{ a(.):[0,\infty) \rightarrow \R^{m_{2}}\,|\, b(.) \mbox{ meas. and } b(t) \in B \mbox{ a.e. }  \}\,.
\end{eqnarray*}
 For  $a(.)\in {\mathcal{A}}$ and $b(.)\in {\mathcal{B}}$  the pay-off  is
\begin{multline*}
 J(x_{0}, a(.), b(.))= \\
 \int_{0}^{\tau}e^{- \lambda t}\;l(x(t; 0, x_{0}, a(.),b(.)),a(t), b(t))dt\,,
\end{multline*}
in which  $\lambda  \geq 0$   (the discount factor) is a given number and $l(.,.,.): \R^{n}\times \R^{m_{1}}\times  \R^{m_{2}}\rightarrow \R$ (the  payoff integrand) is a given function. Here, $\tau$ is the first entry time for $x(t; x_{0},a(.), b(.))$.
\ \\

\noindent
Following Elliot-Kalton \cite{EK72}, we interpret `closed loop' policies for the $a$-player and $b$-player respectively as 
\begin{eqnarray*}
&& \Phi\,:=\, \{\phi: {\mathcal{B}} \rightarrow {\mathcal{A}}\,|\, \phi \mbox{ is non-anticipative}\} ,
\\
&& \Psi\,:=\, \{\psi: {\mathcal{A}} \rightarrow {\mathcal{B}}\,|\, \psi \mbox{ is non-anticipative}\}\, .
\end{eqnarray*}
Here, `$\phi(.)$ is non-anticipative' in the first relation means, `for any $t' \geq 0$,
and $b_{1}(.), b_{2}(.) \in {\mathcal{B}}$,
\begin{multline*}
b_{1}(t)=b_{2}(t)  \mbox{ a.e. } t \in [0,t'] \;\implies \\
\;\phi(b_{1}(.))(t)=\phi(b_{2}(.))(t) \mbox{ for  a.e. }t \in [0,t']\,.
\end{multline*}
`$\psi(.)$ is non-anticipative'
in the second defining relation is analogously defined. Using these interpretations, we define the upper and lower values $u(x)$ and $v(x)$ of the game, for a given starting start $x \in \R^{n}$, to be 
\begin{eqnarray*}
u(x)&=& \underset{\phi \in \Phi}{\sup}\;  \underset{b \in {\mathcal{B}}}{\inf}\; J(x(.;x,\phi (b(.)),, b(.) ),
\\
v(x)&=& \underset{\psi \in \Psi}{\inf}\;  \underset{a \in {\mathcal{A}}}{\sup}\; J(x(.;x, a(.), \psi (a(.))) \;.
\end{eqnarray*}
Define the real valued functions $F(.,.,.)$ and $G(.,.,.)$, with domains in $\R^{n} \times \R \times \R^{n} \rightarrow 
\R$
\begin{eqnarray*}
F(x,u,p)&=& \lambda u + \underset{a \in A}{\inf}\; \underset{b \in B}{\sup}  \;\{p \cdot (-f(x,a,b) - l(x,a,b))\}\,,
\\
G(x,u,p)&=& \lambda u +  \underset{b \in B}{\sup}\;\underset{a \in A}{\inf} \;\{p \cdot (-f(x,a,b) - l(x,a,b))\}.
\end{eqnarray*}
There is an extensive literature on precise conditions on the data, target, etc., under which $u(.)$ coincides with $v(.)$, when $u(.)$ can be characterized as the unique continuous viscosity solution of the HJI (Hamilton Jacobi Isaacs) equation:
\begin{equation}
\left\{
\begin{array}{ll}
\label{HJE_upper}
F(x,u,Du)\,=\,0& \mbox{ for } x \in \R^{n} \backslash {\mathcal{T}},
\\
u(x)=0 \; &\mbox{ for } x \in {\mathcal{T}}\,,
\end{array}
\right.
\end{equation}
and when maximizing closed loop policies for the $a$-player can be obtained from knowledge of $u(.)$. See  \cite{BCD97}, \cite{B98}, \cite{KS88} for expository material on these topics, and \cite{BRP99} for  numerical aspects.\par
\noindent
In this paper, attention focuses on the upper value functon $u(.)$ 
and the associated HJI  equation  (\ref{HJE_upper}). We consider situations in which the target ${\mathcal{T}}$ can be represented as the union of a finite number of closed sets ${\mathcal{T}}_{j}$, $j=1,\ldots,m$:
$$
{\mathcal{T}} \,=\,  \cup_{j=1}^{m} {\mathcal{T}}_{j}\;.
$$
Here, the $b$-player, responding to the closed loop policy of the $a$-player, has a choice over which component ${\mathcal{T}}_{j}$, $j=1,\ldots m$, to exit into, to minimize the payoff. 
Consider the family of `reduced'  value functions $u_{j}(.)$, $j=1,\ldots,m$, that result when the target ${\mathcal{T}}$ is replaced by  the subset ${\mathcal{T}}_{i}$. \par
\noindent
Of interest  are cases in which the value functions $u_{j}(.)$, $j=1,\ldots,m$, for the  target subsets  are
easier to calculate than the value function $u(.)$ for the full target ${\mathcal{T}}$ and when $u(.)$ can constructed as the lower envelope of the $u_{j}(.)$'s, thus:
\begin{equation}
\label{decomp}
u(x)\,= \, \min\{ u_{j}(x)\,|\, j =1,\ldots,m \}.
\end{equation}
The motivation for seeking a decomposition of this nature is as follows.
Optimal control problems are special cases of differential games in which the constraint set $A$ for the $a$-player is a single point; there is then only one possible  open loop policy for the $a$-player, which can therefore be effectively ignored. For optimal control problems, the decomposition (\ref{decomp}) is always valid, since replacing ${\mathcal{T}}$ by one particular ${\mathcal{T}}_{j}$ amounts to a strengthening of the problem constraints, and cannot therefore reduce the value. So, for any $x$ and any $j$, $u(x)\leq u_{j}(x)$. On the other hand, an optimal policy, for the given initial state $x$, must result in the state trajectory exiting into ${\mathcal{T}}_{\bar{k}}$ for some  $\bar k$. But then $u(x) \geq u_{j}(x)$.  These inequalities validate the decomposition (\ref{decomp}).\par
\noindent
When the presence of the $a$-player is restored and we are dealing with a true differential game, decomposition is a much more complicated issue. There are nontheless interesting cases when the decomposition can be achieved. The goal of this paper is to give criteria for decomposition, and to illustrate their application.\par
\noindent
We shall assume that the value functions involved are unique viscosity solutions of the HJI equation with appropriate boundary conditions. This means that checking the validity of the decomposition reduces to answering the question: when is the lower envelope  of a family of viscosity solutions to a particular HJI equation also a viscosity solution? In Section 2 we  give two criteria ($(E)$ and $(C)$)  under which the answer is affirmative. $(E)$ is more general, but $(C)$ is often easier to verify.
$(C)$  is satisfied, in particular, when $F(x,u,.)$ is convex. This is a well-known fact:  the viscosity solution property is preserved under the operation of taking the lower envelopes, for convex Hamiltonians. Notice that, for optimal control problems $F(x,u,.)$ is always convex, so this fact is consistent with the earlier observation that, for optimal control problems, regarded as special cases of differential games, the decomposition is possible. However $(C)$ is weaker than `full' convexity of $F(x,u,.)$, because it requires us to check, for each $x \in \R^{n} \backslash {\mathcal{T}}$,  the convexity inequality only w.r.t.  gradient vectors of the minimizing $u_{j}(.)$'s at $x$. In the examples, this (restricted sense) convexity condition is satisfied while full convexity fails.
%
We provide examples from pursuit/evasion games
in which the decomposition simplifies computations by reducing the state dimensionality.
%
\par
Some examples of the decomposition, without detailed accompanying analysis were presented in \cite{CDCpaper}.
\section{Properties of the Lower Envelope of a Family of Viscosity Solutions}
\label{section2}
Take a function $F(.,.,.): \R^{n}\times \R \times \R^{n} \rightarrow \R$ and consider the partial differential equation
\begin{equation}
\label{1.1}
F(x, u(x), Du(x))\;=\;0\,.
\end{equation}
\begin{definition}
Take an open subset  $ \Omega \subset \R^{n}$   and a  function $u(.):\Omega \rightarrow \R$. Then $u(.)$ is   a continuous viscosity subsolution of (\ref{1.1}) on $\Omega$
if it is continuous and, for each $x\in\Omega$,
\begin{equation}
F(x,u(x), p) \leq 0, \quad \forall p \in D^+u(x)\,. 
\end{equation}
$u(.)$ is   a continuous viscosity supersolution of (\ref{1.1}) on $\Omega$ if it is continuous and, for each
$x\in\Omega$, 
\begin{equation}
F(x,u(x), p) \geq 0, \quad \forall p \in D^-u(x). 
\end{equation}
 $u(.)$   is a continuous viscosity solution of (\ref{1.1}) on $\Omega$ if it is both a continuous subsolution and supersolution of (\ref{1.1}) on $\Omega$.
\end{definition}
Here, $D^{+}u(x)$ and $D^{-}u(x)$ denote, respectively, the  Fr\'echet superdifferential and subdifferential of the continuous function $u(.)$ defined on an open subset of $\R^{n}$ containing the point $x$:
\begin{multline*} D^+u(x):=\\
\left\{p\in\R^N : \limsup_{y\rightarrow x}\frac{u(y)-u(x)-p\cdot (y-x)}{|x-y|} \leq 0\right\}\, ,
\end{multline*}
\begin{multline*}
 D^-u(x):=\\
\left\{p\in\R^N : \liminf_{y\rightarrow x}\frac{u(y)-u(x)-p\cdot (y-x)}{|x-y|}\geq 0\right\}\,.
\end{multline*}
(For the analysis of this paper it is helpful to define continuous viscosity solutions in terms of one-sided Fr\'echet differentials which is equivalent to  the standard definition in terms of gradients of smooth majorizing and minoring functions  \cite{CS04}.
%
\ \\

\noindent
The following proposition gives conditions under which the lower envelope of a collection of continuous viscosity solutions of (\ref{1.1}) is also a continuous viscosity solution, expressed in terms of the limiting superdifferential $\partial^{L}(x)$ of the continuous function $u(.)$ at $x$:
\begin{multline*}
\partial^{L}u(x)\,:=\,\{p\,|\, \exists \mbox{ sequences } \;p_{i}\rightarrow p \mbox{ and }\\ x_{i}\rightarrow x \mbox{ s.t. }
 p_{i}\in D^{+}u(x_{i}) \mbox{ for each } i  \}\,.
\end{multline*}

\begin{proposition}
\label{prop2.2}
Take a collection of  closed sets ${\mathcal{T}}_{j} \subset \R^{n}$, $ j=1,\ldots,m$. For each $j$, let $u_{j}(.)$ be a scalar valued function with domain $\R^{n} \backslash  {\mathcal{T}}_{j}$.
Define 
\begin{multline*}
 I(x)\,=\, \{ j \in \{1,\ldots,m\}\,|\,  u_{j}(x)= \underset{j'}{\mbox{Min}}\;u_{j'} (x)  \} 
\quad \\ \mbox{for each } x \in \R^{n} \backslash ( \cup_{j=1}^m\, {\mathcal{T}}_{j} ) 
\end{multline*}
and 
$$
\Sigma\,=\, \{ x \in \R^{n} \backslash ( \cup_{j=1}^m\, {\mathcal{T}}_{j} ) \,|\,  \mbox{Cardinality}\{I(x)\} \,>\, 1 \} \,.
$$
Take $\bar u(.): \R^{n} \backslash ( \cup_{j=1}^m\, {\mathcal{T}}_{j})\rightarrow \R $ to be the lower envelope function 
\begin{equation*}
\bar u(x)\,=\, \underset{j}{\mbox{Min}}\,\{u_{j}(x) \}\,.
\end{equation*}
\begin{itemize}
\item[(a):] Suppose that $u_{j}(.)$ is a continuous viscosity supersolution of (\ref{1.1}) on $\R^{n} \backslash {\mathcal{T}}_{j} $ for each $j$. Then $\bar u(.)$ is a continuous viscosity supersolution of (\ref{1.1}) on $\R^{n} \backslash ( \cup_{j=1}^m\, {\mathcal{T}}_{j} )$. 
\item[(b):] Suppose that $u_{j}(.)$ is a continuous viscosity subsolution of (\ref{1.1}) on $\R^{n} \backslash {\mathcal{T}}_{j} $ for each $j$,  that $H(.,.,.)$ is continuous and that, for each $x \in \Sigma$, $u_{j}(.)$ is Lipschitz continuous on a neighbourhood of $x$. \par
Consider the hypotheses:
\begin{itemize}
\item[(C):] for any $x \in \Sigma$, any set of vectors $\{p_{j}\,|\, j \in I(x) \}$ such that $p_{j}\in \partial^{L}u_{j}(x)$ for each $j \in I(x)$, and any convex combination $\{ \lambda_{j}\,|\, j \in I(x)\}$,
\begin{eqnarray*}
&& F(x,\bar u (x), \sum_{j \in I(x)} \lambda_{j}p_{j})\,\leq \,
\sum_{j \in I(x)} \lambda_{j} F(x,u_{j}(x),p_{j})\,.
\end{eqnarray*}
\item[(E):] for any $x \in \Sigma$, any set of vectors $\{p_{j}\,|\, j \in I(x) \}$ such that $p_{j}\in \partial^{L}u_{j}(x)$ for each $j \in I(x)$, and any convex combination $\{ \lambda_{j}\,|\, j \in I(x)\}$,
\begin{eqnarray*}
&& F(x,\bar u (x), \sum_{j \in I(x)} \lambda_{j}p_{j})\,\leq \,0
\,.
\end{eqnarray*}
\end{itemize}
\begin{itemize}
\item[(i):]\hspace{0.155 in} $(E)\, \implies \hspace{0.01 in} 
$  `$\,\bar u(.)$ is a continuous viscosity subsolution to (\ref{1.1}) on  $\R^{n} \backslash ( \cup_{j=1}^m\, {\mathcal{T}}_{j} )$'.
\item[(ii):] If, additionally, $u_{j}(.)$ is  $C^{1}$ on a neighborhood of $x$ for each $j$, then\par 
$\,\bar u(.)$ is a continuous viscosity subsolution to (\ref{1.1}) on  $\R^{n} \backslash ( \cup_{j=1}^m\, {\mathcal{T}}_{j} )$  $\implies  (E)$.
\item[(iii):]  
$(C) \,\implies\, (E)$\'.
\end{itemize}
\end{itemize}
\end{proposition}
{\bf Comments.}\par
(i): The proof of the proposition is based on a well-known estimate for  one-sided differentials to lower envelope functions, in terms of the  one-sided differentials to the constituent functions (the `Max Rule').  Such estimates are studied in depth in \cite{LT12}. \\
(ii): The proposition treats separately the preservation of the supersolution and subsolution properties  of viscosity solutions under the operation of taking the lower envelope, because much weaker hypotheses need be imposed in connection wth supersolutions. \\
(iii): We give two sufficient conditions for the lower envelope of a famility of continuous viscosity solutions also to be a continuous viscosity solution, namely $(E)$ and $(C)$. 
$(C)$ is a more restrictive condition, but it is useful because, as illustrated in the following examples, it can be easier to verify. \\
(iv): The proposition is an analytical tool for decomposing a differential game (associated with the value function $\bar u(.)$) into a collection of simpler problems. The critical hypothesis in this proposition is $(E)$ (or $(C)$). $(C)$ is automatically satisfied when $F(x,u,.)$ is convex. This special case of the proposition is well-known \cite{CS04}. 
However $(C)$ imposes a convexity type condition on $F(x,u,.)$, {\it only with respect to selected vectors in its domain}. In some cases, examples of which given below, the restricted sense convexity hypothesis is satisfied but the full convexity hypothesis is violated; the proposition thereby identifies a new class of differential games for which the decomposition is possible.\par
\noindent
{\bf Proof of Prop. \ref{prop2.2}}.

\noindent
(a): Suppose that $u_{j}(.)$ is a continuous viscosity supersolution of (\ref{1.1}) on $\R^{n} \backslash {\mathcal{T}}_{j} $ for each $j$. Take any $x \in \R^{n} \backslash ( \cup_{j=1}^m\, {\mathcal{T}}_{j} )$ and $p \in D^{-}\bar u(x)$. Then
$$
\bar u(x')  - \bar u(x) \,\geq\, p\cdot (x' -x) -\mbox{{\it o}}(|x' -x|)\, ,
$$
for all $x'  \in \R^{n} \backslash ( \cup_{j=1}^m\, {\mathcal{T}}_{j} )$.
(Here,  $o(.):\R^{+} \rightarrow \R^{+}$ is some function such that $\lim_{s \downarrow 0} o(s)/s \rightarrow 0$.) Choose any $j \in I(x)$. We know that $u_{j}(x)= \bar u(x)$ and $u_{j}(x') \geq \bar u(x')$ . It follows that, for all  $x' \in \R^{n} \backslash ( \cup_{j=1}^m\, {\mathcal{T}}_{j} )$, 
$$
u_{j}(x')  - u_{j}(x) \,\geq\, p\cdot (x' -x) -\mbox{{\it o}}(|x' -x|)\,.
$$
But then $p \in D^{-}u_{j}(x)$ and, since $u_{j}$ is a continuous viscosity supersolution, we have $F(x, u_{j}(x), p)\geq 0$. It follows that $F(x, \bar u(x), p)\geq 0$. Since $\bar u(. )$ is continuous, we have established that $\bar u(.)$ is a continuous viscosity subsolution of (\ref{1.1}) on  
$\R^{n} \backslash ( \cup_{j=1}^m\, {\mathcal{T}}_{j} )$.
\ \\

\noindent
(b)(i): Suppose that $u_{j}(.)$ is a continuous viscosity subsolution of (\ref{1.1}) on $\R^{n} \backslash {\mathcal{T}}_{j} $ for each $j$. Take any $x \in \R^{n} \backslash ( \cup_{j=1}^m\, {\mathcal{T}}_{j} )$ and $p \in D^{+}\bar u(x)$. We must show that
\begin{equation}
\label{1.2}
F(x,\bar u(x),p ) \,\leq \,0\,.
\end{equation}
Suppose first that $x\notin \Sigma$, i.e. $I(x)$ contains a single index value $j$. Then, since the $u_{i}(.)$'s are continuous, $\bar u(x')= u_{j}(x')$ for all $x'$ in some neighbourhood of $x$. It follows that $p\in D^{+}u_{j}$ and so $F(x,u_{j}(x),p )(=F(x,\bar u (x),p )) \,\leq \,0$. We have confirmed (\ref{1.2}) in this case.\par
\noindent
It may be assumed then that $x \in \Sigma$. Now, $u_{j}(.)$ is Lipschitz continuous on a neighbourhood of $x$ for each $j \in I(x)$. Since $p\in D^{+}\bar u(x)$, it is certainly the case that $p \in \partial^{L} \bar{u}(x)$. Using the property that $\bar u(x')$ coincides with $\max \{ u_{j}(x')\,|\, j \in I(x')   \}$ for $x'$ in some neighbourhood of $x$, we deduce from the Max Rule for limiting subdifferentials of Lipschitz continuous functions  (see, e.g., \cite[Thm. 5.5.2]{V00}) applied to $-\bar u(.)$ the following representation for $p$:
$$
p =\sum_{j\in I(x)}\lambda_{j}p_{j}\,,
$$
for some convex combination $\{\lambda_{j}\,|\, j \in I(x)\}$ and vectors $p_{j} \in \partial^{L}u_{j}(x)$, $j\in I(x)$. But then, by hypothesis (E),
$$
F(x,\bar u(x),p)= F(x,\bar u(x),\sum_{j \in I(x)}\lambda_{j} p_{j}) 
\,\leq 0\,.
$$
We have confirmed (\ref{1.2}) and so (b)(i) is true.\par
\noindent
(b)(ii): Take any $x \in \Sigma$. Suppose that the $u_{j}$'s are continously differentiable of a neighbourhood of $x$ and that $\bar{u}(.)$ is a viscosity solution. Take any convex combination $\{\lambda_{i}\}$ on $I(x)$. Then, for all $x'$ in some neighborhood of $x$,
\begin{multline*}
\bar{u}(x') - \bar{u}(x) \leq \sum_{i \in I(x)} \lambda_{i} (u_{i}(x') -u_{i}(x)) \\
\leq
  \sum_{i\in I(x)} \lambda_{i}\nabla u_{i}(x) \cdot (x' -x) +o(|x'-x|)\;.
\end{multline*}
This last inequality tells us that $\sum_{i\in I(x)} \lambda_{i}\nabla u_{i}(x)$ is a limiting superdifferential of $\bar{u}(.)$ at $x$. But then, since $\bar{u}(.)$ is a viscosity subsolution,
$$
F(x,\bar u(x),\sum_{j \in I(x)}\lambda_{j} p_{j}) 
\,\leq 0\,.
$$
We have confirmed that $(E)$ is true.\par
\noindent 
(b)(iii): Take any convex combination $\{\lambda_{i}\}$ on $I(x)$ and vectors $p_{i}\in \partial^{L}u_{i}(x) $ for $i \in I(x)$. 
It follows from the definition of the limiting supergradient that, for each $i$, there exist sequences $x_{j}^{i} \rightarrow x$ and $p_{j}^{i}\rightarrow p_{i}$ such that $p_{j}^{i} \in D^{+}u_{i}(x_{j}^{i})$ for $i=1,2, \ldots$ But then, for each $i\in I(x)$,
$$
F(x_{j}^{i}, u^{j}(x_{j}^{i}),  p_{j}^{i}) \leq 0, \; i=1,2,\ldots,
$$
since  the $u_{i}(.)$'s are viscosity subsolutions.
It follows that $\sum_{j \in I(x)}\lambda_{j} F(x_{j}^{i}, u^{j}(x_{j}^{i}),  p_{j}^{i}) \leq 0, \; i=1,2,\ldots$. Noting the continuity of $F(.,,.,)$ and also the $u_{j}(.)$'s, we may pass to the limit  as $i \rightarrow \infty$ to obtain
$$
\sum_{j \in I(x)}\lambda_{j} F(x,u_{j}(x), p_{j})\,\leq 0\,.
$$
Assume $(C)$. Then
$$
F(x,\bar u(x),\sum_{j \in I(x)}\lambda_{j} p_{j}) \leq \sum_{j \in I(x)}\lambda_{j} F(x,u_{j}, p_{j})\,\leq 0\,,
$$
which is $(E)$. 
\par
\noindent
\section{Pursuit Evasion Games} 
\label{section3}
Pursuer/evader games are examples of the game posed in the Introduction.  There is an extensive literature on such games, going back to Rufus Isaacs' work in the 1960's, and his monograph  \cite{I65} contains many examples. Expository material is to be found in    \cite{F71},  \cite{KS88}.  We note also  \cite{C89}, \cite{IL81}, \cite{I05}, \cite{IL81}, \cite{P76}, and \cite{VSH02}. But none of these references systematically address decomposions of the game, each element of which is generated by a target subset. Pursuer/evader games is an application area for the methods proposed in this paper; they provide exemplar problems, both where decomposion is possible, and where it is not.
\par
\noindent
We consider  zero sum differential games which terminate when one of the pursuers is sufficently close to one of the evaders, where `closeness' is understood in the sense of a specified target. The pay-off is the time until the target is attained. We analyse a number of examples, involving different numbers of pursuers and evaders, and different targets.\par
\noindent
 The $a$-player is the collection of $m_{1}$ evaders, labelled $1,\ldots, m_{1}$, and the $b$-player the collection of $m_{2}$ pursuers, labelled  $m_{1}+1, \ldots, m_{1}+m_{2}$.
The states of individual pursuers and evaders $x_{1},\ldots, x_{m_{1}}$ and $x_{m+1}\ldots, x_{m_{1}+m_{2}}$ are governed by the equations
\begin{eqnarray*}
&& \frac{dx_{1}}{dt} \;=\; f_{1}(x_{1},a_{1})\,, \,.\, .\, .\,  \,, \,
\frac{dx_{m_{1}}}{dt} \;=\; f_{m_{1}}(x_{m_{1}},a_{m_{1}}) )
\\
&&\frac{dx_{m_{1}+1}}{dt} \;=\; f_{m_{1}+1}(x_{m_{1}+1},b_{1}),\; \ldots \;,\\
&& \frac{dx_{m_{1}+m_{2}}}{dt} \;=\; f_{m_{1}+m_{2}}(x_{m_{1}+m_{2}},b_{m_{2}})\;.
\end{eqnarray*}  
The variables $a_{1},\ldots,a_{m_{1}}$ and $b_{1},\ldots,,b_{m_{2}}$ are interpreted as controls for the evaders and the pursuers, respectively,  which are subject to the constraints
$$
a_{i}\in A_{i},\; i= 1,\ldots ,m_{1} , \mbox{ and }  b_{i} \in B_{i},\;  i= 1,\ldots ,m_{2}\;.
$$
Here, $f_{i}(.,.): \R^{n} \times \R^{r_{i}}\rightarrow \R^{n}$, $i=1,\ldots,m_{1}+m_{2}$ are given functions,  and $A_{i} \subset \R^{r_{i}}$, $1,\ldots,m_{1}$ and $B_{i} \subset \R^{r_{i+m_{1}}}$, $1,\ldots,m_{2}$, are given subsets.\\
\noindent
We regard $a_{1},\ldots, a_{m_{1}}$ and $b_{1},\ldots, b_{m_{2}}$ as block components of a single evader  control and pursuer control respectively. Take the state to be $x=\mbox{col}\, \{x_{1}\ldots, x_{m_{1}+m_{2}}\}$. The open loop policy spaces for evader and pursuer are
\begin{eqnarray*}
 {\mathcal{A}}\,&:=&\, \{\mbox{meas. mappings }a_{i}: [0,\infty) \rightarrow \R^{r_{i}}, \, 
i=1,\ldots,m_{1}\,|\\
\, & &a_{i}(t)\in A_{i} \mbox{ a.e. } \mbox{for each }i\},
\\
 {\mathcal{B}}\,&:=&\, \{\mbox{meas. mappings }\;b_{i}: [0,\infty) \rightarrow \R^{r_{i+m_{1}}}, \, 
i=i,\ldots,  m_{2}\,|\\
\, & & b_{i}(t) \in B_{i} \mbox{ a.e. for each }i  \}\,.
\end{eqnarray*}
Write $\Phi$ for the space of non-anticipative mappings $\phi: {\mathcal{B}} \rightarrow {\mathcal{A}}$. The game fits the formulation  Section 1, with $\lambda =0$, and may be summarized as:
$$
(P')
\left\{
\begin{array}{l}
\underset{\phi \in \Phi}{\mbox{Maximize}} \;
\underset{\{b_{i}\} \in {\mathcal{A}}}{\mbox{Minimize}} \int_{0}^{\tau} \,1 \,dt
\\
\left\{
\begin{array}{l}\dot x_{1}(t)\;=\; f_{1}(x_{1}(t),a_{1}(t))\\
...\\
\dot{x}_{m_{1}+ m_{2}}(t)= f_{m_{1}+m_{2}}(x_{m_{1}+m_{2}},b_{m_{2}}),\, \mbox{ a.e.}

\end{array}\right. \\
(a_{1}(t),...,a_{m_{1}},b_{1}(t),...,b_{m_{2}}(t))\\
\phantom{bhbrgfhr4he}\in A_{1}\times ... \times  A_{m_{1}} \times B_{1}\times ... \times  B_{m_{2}},
\mbox{ a.e.}
\\
\mbox{in which } (a_{1}(.),\ldots, a_{m_{1}}(.))= \phi(b_{1}(.),\ldots, b_{m_{2}}(.))\\
 \mbox{ and }
\tau \mbox{ is first entry time into } {\mathcal{T}}\\
(x_{1}(0),\ldots, x_{m_{1}+m_{2}}(0))\,=\, (\bar x_{1},\ldots,\bar x_{m_{1}+ m_{2} })
\end{array}
\right.
$$
for some given $(\bar x_{1},\ldots,\bar x_{m_{1}+m_{2}}) \in \R^{n}\times \ldots \times \R^{n}$. Here ${\mathcal{T}}$ is a given closed subset of $\R^{n}\times \ldots \times \R^{n} $.
The Hamilton-Jacobi-Isaacs equation is
\begin{equation}
\label{HJE}
F(x_{1},\ldots,x_{m_{1}+m_{2}}, D_{x_{1}}u,\ldots, D_{x_{m_{1}}+m_{2}} u )\,=\,0\, ,
\end{equation}
in which 
\begin{multline*}
F(x_{1},\ldots,x_{m_{1}+m_{2}}, p_{1},\ldots, p_{m_{1}+m_{2}} )\,=\\
\, 
- \sum_{i=1}^{m_{1}}H^{i}(x_{i},p_{i}) \,+\,\sum_{i=m_{1}+1}^{m_{1}+m_{2}} H^{i}(x_{i},-p_{i})  -1.
\end{multline*}
\begin{multline}
\label{small}
\mbox{Here }\;H^{i}(x_{i},p_{i}):= \\
\left\{
\begin{array}{ll}
\underset{a_{i} \in A_{i}}{\sup}\, p_{i}\cdot f(x_{i},a_{i})&\mbox{for } i =1,\ldots m_{1}
\\
\underset{b_{i-m_{1}} \in B_{i-m_{1}}}{\sup}\,
p_{i}\cdot f(x_{i},b_{i-m_{1}})&\mbox{for } i =m_{1}+1,\ldots m_{1}+ m_{2}.
\end{array}
\right.
\end{multline}
\subsection{A Single Pursuer/Multiple Evaders Game}
Consider first a case of the pursuit/evasion game, written $(P^{1})$, in which $m_{1}=m >1$, $m_{2}=1$ and $n=1$ (a single pursuer/multiple evaders game in 1D space). The states of  the $m$ evaders, labeled $1,\ldots, m$ and of the one pursuer, labeled $m+1$, are interpreted as the positions of the evaders and pursuer. The game terminates when the pursuer is first at a distance $r$ from one of the evaders, where $r \geq0$ is a given constant. Accordingly, we take
$$
{\mathcal{T}} \,=\, {\mathcal{T}}_{1} \cup \ldots \cup{\mathcal{T}}_{m}\, ,
$$
in which, for $i =1,\ldots,m$,
$$
{\mathcal{T}}_{i}\,:=\, \{ (x_{1}, \ldots, x_{m+1}) \,|\,     |x_{m+1}-x_{i}| \leq r   \}\;.
$$
The Hamilton-Jacobi-Isaacs equation is
\begin{equation}
\label{HJE1}
F^{1}(x_{1},\ldots,x_{m+1},D_{x_{1}}u,\ldots, D_{x_{m+1}}u)\,=\,0\, ,
\end{equation}
in which
\begin{multline*}
F^{1}(x_{1},\ldots,x_{m+1}, p_{1},\ldots, p_{m+1} )\,=\\
\, 
- \sum_{i=1}^{m}H^{i}(x_{i},p_{i}) \,+\, H^{m+1}(x_{m+1},-p_{m+1})  -1\,,
\end{multline*}
where
\begin{equation}
\label{small'}
\begin{split}
H^{i}(x_{i},p_{i})\,:=\,\underset{a_{i} \in A_{i}}{\sup}\, p_{i}\cdot f(x_{i},a_{i})\; i =1,\ldots m,\,\\
H^{m+1}(x_{m+1},p_{m+1})\,=\,
\underset{b_{1} \in B_{1}}{\sup}\,
p_{m+1}\cdot f(x_{m+1},b_{1})\,.
\end{split}
\end{equation}
Now take $(P^{1}_{i})$ to be the modification of  $(P^{1})$, when  ${\mathcal{T}}_{i}$ replaces ${\mathcal{T}}$, $i=1,\ldots,m$. 
Let us assume that, for each $i$, the value function $u_{i}(.)$ for $(P^{1}_{i})$ is a continuous viscosity solution of (\ref{HJE}). The following proposition tells us that we can construct a viscosity solution  to (\ref{HJE1}) from the $u_{i}(.)$'s, by taking the pointwise infimum.
\begin{proposition}
\label{prop_3_1}
For  $ i=1,\ldots m$, let $u_{i}(.)$ be the upper value for $(P^{1}_{i})$.
Assume
\begin{itemize}
\item[(a):] For  $i=1,\ldots m$, $u_{i}(.)$ is a continuous viscosity solution of $(\ref{HJE})$ on $\R^{m+1}\backslash {\mathcal{T}}_{i}$.\smallskip
\item[(b):]For any $i,j \in \{1,\ldots,m \}$, $i \not=j$, and $(x_{1}, \ldots,x_{m+1}) \in \R^{m+1} \backslash {\mathcal{T}}$ such that $u_{i}(x_{1},\ldots,x_{m+1})= u_{j}(x_{1},\ldots,x_{m+1})$, $u_{i}(.)$ and $u_{j}(.)$ are Lipschitz continuous on a neighborhood of $(x_{1}, \ldots,x_{m+1})$.
\end{itemize}
Then
\begin{multline*}
\bar u(x_{1},\ldots,x_{m+1}) := \\
\min \{ u_{1}(x_{1},\ldots x_{m}),\ldots,  u_{m+1}(x_{1},\ldots,x_{m+1})   \}
\end{multline*}
is a continuous viscosity solution of (\ref{HJE}) on  $(\R\times \ldots \R)\backslash {\mathcal{T}}$.
\end{proposition}
\vspace{0.05 in}

\noindent
{\bf Comment.}
{\it Suppose hypotheses are imposed, ensuring that (1): for each $i$, the HJI equation for $(P^{1}_{i} )$ has a continuous viscosity solution on $(\R^{n}\times \ldots \R^{n})\backslash {\mathcal{T}}_{i}$ with a continous extension to ${\mathcal{T}}_{i}$, on which set the solution vanishes, and (2):  the value function  $(P^{1})$ is the unique continuous viscosity solution on $(\R^{n}\times \ldots \R^{n})\backslash {\mathcal{T}}$ that has a continous 
extension to ${\mathcal{T}}_{i}$, on which set the solution vanishes. The proposition tells us that, under these circumstances, the upper value $u(.)$ for $(P^{1})$ can be calculated as the lower envelope of the continuous viscosiy solutions for the $(P^{1}_{i})$'s. (Notice that, since all upper values concerned are non-negative, and each $u_{i}(.)$ is assumed to have a continuous extension to ${\mathcal{T}}_{i}$, on which set it vanishes, the lower envelope has a continuous extension to ${\mathcal{T}}$, on which set it vanishes.)}
\ \\

\noindent
{\bf Proof of Prop. \ref{prop_3_1}.} Note that, for any $i$,   $u_{i}(x_{1}, ... , x_{m+1})$ depends only on the two variables $(x_{i}, x_{m+1})$. This is because the first entry time into ${\mathcal{T}}_{i}$ only concerns the state trajectories associated the $i$'th evader and the pursuer (labelled $m+1$). \par
\noindent 
In view of the hypotheses imposed on the $u_{i}(.)$'s, the fact that $\bar u(.)$ is a viscosity solution of (\ref{HJE}) will follow from Prop. \ref{prop2.2}, if we can confirm hypothesis (C) of this proposition. Take any $z =(x_{1}, \ldots,x_{m+1})\in \R^{m+1} \backslash {\mathcal{T}}$, any index set $I(z)$ (of cardinality $l >1$) such that the values $u_{i}(z)$, $i\in I(z)$, coincide, and any convex combination $\{\lambda_{i}\}$ from $I(z)$. To simplify, assume index values have been re-ordered so that $I(z)= \{1,\dots,,l \}$.
Take also  $\tilde p^{i}\in \R^{n} $, $i =  1,\ldots,l$  such that
\begin{equation}
\label{special_a1}
\tilde p^{i} \,:=\, (0,\ldots,0,p^{i}_{i},0,\ldots,0, p^{i}_{m+1}) \in  \partial^{L}u_{i}(z)\,.
\end{equation}
(The possibly non-zero components  $p^{i}_{i}$  and $p^{i}_{m+1}$ of $\tilde p^{i}$ appear at the  $i$'th and $(m+1)$'th locations. We must show $\eta(\lambda_{1}, \ldots,\lambda_{l}) \geq 0$, where
\begin{eqnarray*}
&& \eta(\lambda_{1}, \ldots,\lambda_{l})\,:=\,
\sum_{i=1}^{l}\lambda_{i}F(z, \tilde p_{i})- F(z, \sum_{i=1}^{l}\lambda_{i }\tilde p_{i})\,.
\end{eqnarray*}
Noting the special structure  (\ref{special_a1}) of the $\tilde p_{i}$'s and the fact that $H^{i}(x_{i},p_{i})=0$ when $p_{1}=0$, for each $i$, we see that
\begin{multline*}
 \eta(\lambda_{1}, \ldots,\lambda_{l})=\,\\
\,
 \sum_{i=1}^{l} \lambda_{i}\left(H^{m+1}(x_{m+1},-p^{i}_{m+1}) - H^{i}(x_{i},p^{i}_{i}) \right)
\\
 -\left( 
H^{m+1}(x_{m+1},-\sum_{i=1}^{l} \lambda_{i}  p^{i}_{m+1}) \,-\,\sum_{i=1}^{l}H^{i}(x_{i}, \lambda_{i}p^{i}_{i}))
\right) \;.
\end{multline*}
We achieve a further simplification from the fact that $H^{i}(x_{i}, .)$ is positively homogeneous, so $H^{i}(x_{i}, \lambda_{i} p^{i}_{m+1})= \lambda_{i }H^{i}(x_{i},p^{i}_{m+1}) $. This gives
\begin{multline*}
 \eta(\lambda_{1}, \ldots,\lambda_{l})\,=\\
\,
 \sum_{i=1}^{l}\ \lambda_{i}H^{m+1}(x_{m+1},-p_{m+1}^{1}) \,-\,
H^{m+1}(x_{m+1},-\sum_{i=1}^{l} \lambda_{i} \ p^{i}_{m+1})
\;.
\end{multline*}
But then $\eta(\lambda_{1}, \ldots,\lambda_{l})$ is non-negative, because the term $H^{m+1}(x_{m+1},.)$, defined by (\ref{small}), is convex.  The proof is complete.

\subsection{A Multiple Pursuers/Single Evader Game} 
\label{section4_2}
Consider next a case of the pursuit/evader game, written $(P^{2})$, in which $m_{1}=1$, $m_{2}\geq 1$ and $n=2$ (single pursuer/multiple evaders).  The dynamic behavior of each player is modelled as a thrust acting on a mass, in 1D space, with saturating damping. The state equations, governing the position and velocity of each player, are taken to be,  for $i=2,\ldots, m+1$,
\begin{small}
$$
\left[
\begin{array}{c}
\dot x^{1}_{1}
\\
\dot x^{1}_{2} 
\end{array}
\right]
= 
\left[
\begin{array}{c}
x^{1}_{2}
\\
-d_{1}(x^{i}_{2}) + a_{1} 
\end{array}
\right]\,\mbox{ and } \;
\left[
\begin{array}{c}
\dot x^{i}_{1}
\\
\dot x^{i}_{2} 
\end{array}
\right]
= 
\left[
\begin{array}{c}
x^{i}_{2}
\\
-d_{i}(x^{i}_{2}) + b_{i-1} 
\end{array}
\right].
$$
\end{small}
Here, $d_{i}(.):\R \rightarrow \R$, $i=1,\ldots, m+1$ are given functions
 satisfying
\begin{equation}
\label{d_prop}
|d_{i}(y)-d_{i}(y')| \leq k_{d}|y-y'|,\; d_{i}(y) \leq c_{d},
\end{equation}
for all $ y,y' \in \R\,\mbox{ and } i=1,\ldots, m+1$ for some constants $k_{d}>0$ and $c_{d}>0$.
The control actions the players are required to satisfy
\begin{equation}
\label{bounds}
|a| \leq \alpha \mbox{ and } |b| \leq \beta_{i}\; \mbox{ for } i=1,\dots. m\,.
\end{equation}
for positive constants $\alpha$, $\beta_{1},\ldots, \beta_{m}$.
We assume that
\begin{equation}
\label{constraints}
\beta_{i}>  \alpha + 2 \times c_{d}\, \mbox{ for i=1,\ldots, m}.
\end{equation}
The game terminates when one of the pursuers overtakes the evader. Thus, we take the target to be 
$$
{\mathcal{T}} \,=\, {\mathcal{T}}_{1} \cup \ldots \cup{\mathcal{T}}_{m}\, ,
$$
in which, for $i =1,\ldots,m$,
\begin{multline*}
\mathcal{T}_{i}\,:=
\, \{ (x^{1}=(x_{1}^{1},x_{2}^{1}), \ldots, \\
x^{m+1}=(x^{m+1}_{1},x^{m+1}_{2}) \,|\,     x^{i}_{1} \geq  x^{1}_{1}  \}\;.
\end{multline*}
The HJI equation is 
\begin{equation}
\label{HJE2}
F^{2}(x^{1},\ldots,x^{m+1},D_{x^{1}}u,\ldots, D_{x^{m+1}}u)\,=\,0\, ,
\end{equation}
in which
\begin{multline*}
F^{2}(x^{1},\ldots,x^{m+1}, p^{1},\ldots, p^{m+1} )\,=\, \\
\left( 
 \sum_{i=1}^{m+1} \left( 
-p^{i}_{1}x^{i}_{2} -p^{i}_{2}d(x^{i}_{2})
\right)
\right)
-\alpha \times |p^{1}_{2}| +  \sum_{i=2}^{m+1} (\beta_{i} \times |p^{i}_{2}|)\, .
\end{multline*}
%
%
Let $(P^{2}_{i})$ to be the modification of  $(P^{2})$, when the target ${\mathcal{T}}_{i}$ replaces ${\mathcal{T}}$, $i=2,\ldots,m+1$. 

\begin{proposition}
Let $u_{i}(.)$ be the upper value for $(P^{2}_{i})$, for $ i=2,\ldots m+1$. 
Assume
\begin{itemize}
\item[(a):] For $i=2,\ldots m+1$, $u_{i}(.)$ is a continuous viscosity solution of $(\ref{HJE})$ on $(\R^{2})^{m+1}\backslash {\mathcal{T}}_{i}$.
\smallskip
\item[(b):]For any $i,j \in \{2,\ldots,m+1 \}$, $i \not=j$, and $(x^{1}, \ldots,x^{m+1}) \in (\R^{2})^{m+1} \backslash {\mathcal{T}}$ such that 
$u_{i}(x^{1},\ldots,x^{m+1})= u_{j}(x^{1},\ldots,x^{m+1})$, $u_{i}(.)$ and $u_{j}(.)$ are Lipschitz continuous on a neighborhood of $(x^{1}, \ldots,x^{m+1})$.
\end{itemize}
Then
\begin{multline*}
\bar u(x^{1},\ldots,x^{m+1}) := \\
\min \{ u_{1}(x^{1},\ldots x^{m+1}),\ldots,  u^{m+1}(x_{1},\ldots,x^{m+1})   \}
\end{multline*}
is a continuous viscosity solution of (\ref{HJE2}) on  $(\R^{2})^{m+1}\backslash {\mathcal{T}}$.
\end{proposition}
\vspace{0.05 in}

\noindent
{\bf Comment.}
{\it When, for each $i$, the HJI equation for $(P^{2}_{i} )$ has a continuous viscosity solution $u_{i}(.)$ on $(\R^{2})^{m+1}\backslash {\mathcal{T}}_{i}$ (with appropriate boundary values) and the value function $u(.)$ for  $(P^{2})$  is the unique continuous viscosity solution on $(\R^{2})^{m+1}\backslash {\mathcal{T}}$ (with appropriate boundary values), the proposition describes how the value function for $(P^{2} )$ can be obtained, as the pointwise infimum of the $u_{i}(.)$'s.}
\ \\

\noindent
{\bf Proof.}  Note that, for $i=2,\ldots, m+1$,   $u_{i}(x^{1}, \ldots, x^{m+1})$ depends only on the two variables $(x^{1}, x^{i})$,  since the first entry time into ${\mathcal{T}}_{i}$ only concerns the state trajectories associated with the $i$'th pursuer and the evader.  We write  $u_{i}(x^{1},x^{i})$, suppressing irrelevant arguments in the notation. Note that by assumptions (\ref{d_prop}) and (\ref{constraints}) (which tell us that all evaders can accelerate at a faster rate than the evader), $u_{i}(x^{1},x^{i})$ is finite when $x^{1}_{1} \geq x^{i}_{1}$.
\ \\

\noindent 
The left side of the HJI equation $F^{2}=0$ can be decomposed as
\begin{equation}
\label{3_HJI}
F^{2} = F^{21} + F^{22}\, ,
\end{equation}
where $F^{21}$ and $F^{22}$, evaluated at $((x^{1}_{1}, x^{1}_{2}),\ldots,$ $(x^{m+1}_{1}, x^{m+1}_{2}),$ $(p^{1}_{1},p^{1}_{2}),\ldots, (p^{m+1}_{1},p^{m+1}_{2}))$, are:
\begin{equation}
\label{F21}
F^{21}\,=\,  - \sum_{i=1}^{m+1} (p^{i}_{1}x^{i}_{2}-p^{i}_{2}d( x^{i}_{2})) \, ,
\end{equation}
and
\begin{equation}
\label{F22}
F^{22} \,= \,
  - \alpha \times |p^{1}_{2}|  + \sum_{i=2}^{m+1} \beta_{i} \times |p^{i}_{2}|\, .
\end{equation}

\noindent
We shall make use of the following Lemma, whose proof appears in the appendix.

\begin{lemma} 
\label{control}
Let $\bar a(.)$ be  the open loop strategy $\bar a(.) \equiv +1$ for the evader, and let $a(.)$ be any other open loop strategy. Take initial states for the evader $z= (z_{1},z_{2})$ and $z' = (z^{'}_{1},z^{'}_{2})$ such that $z_{1}^{'} \geq z_{1}$ and $z_{2}^{'} \geq z_{2}$. Then 
\begin{equation}
\label{monotone}
x^{1}_{1}(t; \bar a(.),z') \geq x^{1}_{1}(t; a(.),z) \quad \mbox{for all } t \geq 0\,,
\end{equation}
where $t \rightarrow (x_{1}^{1},x_{2}^{1})(t; a(.),z)$ is the state trajectory for the evader, under the open loop strategy $a(.)$ and for  initial state $z$.
\end{lemma}
Fix $i$, and consider $(P^{2}_{i})$. We deduce from the lemma  that the optimal closed loop strategy for the $a$-player (the evader) is $\bar \phi(b_{i}(.)) \equiv +1$, for arbitrary initial state $(x^{1}=(x^{1}_{1},x^{1}_{2}),  x^{i}_{1}=(x^{i}_{1},x^{i}_{2}))$ such that $x^{1}_{1} > x^{i}_{1}$. Furthermore, if the $a$-player applies this optimal strategy then, for any open loop strategy $b_{i}(.)$, the effect of increasing the $x^{1}_{2}$ component of the initial state is to increase the first interception time. We conclude that
\begin{equation}
\label{mono}
x_{2}^{1} \rightarrow u_{i}((x^{1}_{1},x^{1}_{2}),  (x^{i}_{1},x^{i}_{2}))
\quad \mbox{is monotone increasing}\;
\end{equation}
for arbitrary $((x^{1}_{1},  (x^{i}_{1},x^{i}_{2})) \in \R^{3}$ such that
$x^{1}_{1} > x^{i}_{1}$.
\ \\

\noindent
Once again, we shall deduce that the lower envelope $\bar u(.)$ of the $u_{i}$'s is a continuous viscosity solution (\ref{HJE2})  from   Prop. \ref{prop2.2}, by verifying hypothesis (C). Take any $z =(x^{1}, \ldots,x^{m+1})\in (\R^{2} \times \ldots \R^{2})\backslash {\mathcal{T}}$, any index set $I(z)$ (of cardinality $l >1$) such that the values $u_{i}(z)$, $i\in I(z)$ coincide, and any convex combination $\{\lambda_{i}\}$ from $I(z)$. We may assume that index values have been re-ordered so that $I(z)= \{2,\dots,,l+1 \}$.
For $i =  2,\ldots,l+1$, take any  $\tilde p^{i}\in \R^{n} $,   such that
\begin{multline}
\label{special}
\tilde p^{i} \,:=\, ((p^{i,1}_{1}, p^{i,1}_{2}), (0,0),\ldots,\\
 (0,0), (p^{i,1}_{1}, p^{i,1}_{2}),(0,0), \ldots (0,0)) \in  \partial_{L}u_{i}(z)\, .
\end{multline}
(We have used the fact that $u_{i}$ depends only on $(x^{1}=(x^{1}_{1}, x^{1}_{2}),x^{i}=(x^{i}_{1},x^{i}_{2}))$.) The possibly non-zero components  $(p^{i,1}_{1},p^{i,1}_{2})$  and $(p^{i}_{1},p^{i}_{2})$ of $\tilde p^{i}$ appear at the first and $i$'th locations. Note that, by (\ref{mono}), 
\begin{equation}
\label{samesign}
p^{i,1}_{2} \,\geq \,0\quad \mbox{for } i=2,\ldots, m+1\,.
\end{equation}
Verification of hypothesis (C) requires us to  show that $\eta(\lambda_{2},\ldots,\lambda_{l+1}) \geq 0$, where
\begin{eqnarray*} \eta(\lambda_{2}, \ldots,\lambda_{l+1})&:=&
\sum_{i=2}^{l+2}\lambda_{i}F^{21}(z, \tilde p_{i})- F^{21}(z, \sum_{i=2}^{l+2}\lambda_{i }\tilde p_{i})
\\
&& \hspace{-0.1 in} \,+\,
\sum_{i=2}^{l+2}\lambda_{i}F^{22}(z, \tilde p_{i})- F^{22}(z, \sum_{i=2}^{l+2}\lambda_{i }\tilde p_{i})\, .
\end{eqnarray*}
Because $F^{21}(z,.)$ is  linear, we have
\begin{eqnarray*} 
\eta(\lambda_{2}, \ldots,\lambda_{l+1})
&:=&
\sum_{i=2}^{l+2} \lambda_{i} F^{22}(z, \tilde p_{i})- 
F^{22} (z, \sum_{i=2}^{l+2}\lambda_{i } \tilde p_{i})
\\
&=& c_{1}+ c_{2}\,,
\end{eqnarray*}
where
\begin{eqnarray*}
c_{1}\,&:=&\,\sum_{i=2}^{l+1}\left( -\alpha \lambda_{i}|p^{i,1}_{2}| + \alpha| \lambda_{i} p^{1}_{2}| \right)
\; \mbox{ and }\;\\
c_{2}\,&:=&\, \sum_{i=2}^{l+1}\lambda_{i}f^{22}(p^{i}) - f^{22}(z,\sum_{i=2}^{l+1}\lambda_{i}p^{i} )\,,
\end{eqnarray*}
in which 
$$
f^{22}(((p^{1}_{1}, p^{1}_{2}),\ldots,(p^{m+1}_{1}, p^{m+1}_{2})\,:=\,   \sum_{i=2}^{m+1} \beta_{i} \times  |p^{i}_{2}|\,.
$$
But $c_{1}=0$ since, by (\ref{samesign}), the $p^{i,1}_{2}$'s all have the same sign. Also, $c_{2} \geq 0$, by convexity of $f^{22}(.)$. We have confirmed $\eta(\lambda_{2},\ldots,\lambda_{l+1}) \geq 0$, and the proof of the proposition is complete.
\par
\noindent
\begin{figure}[t]
\begin{center}
\includegraphics[height=5.6cm]{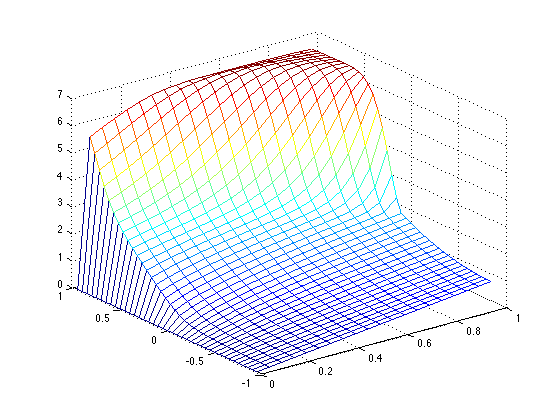}
\caption{Value function for a one-pursuer one-evader game}\label{fig10}
\includegraphics[height=5.6cm]{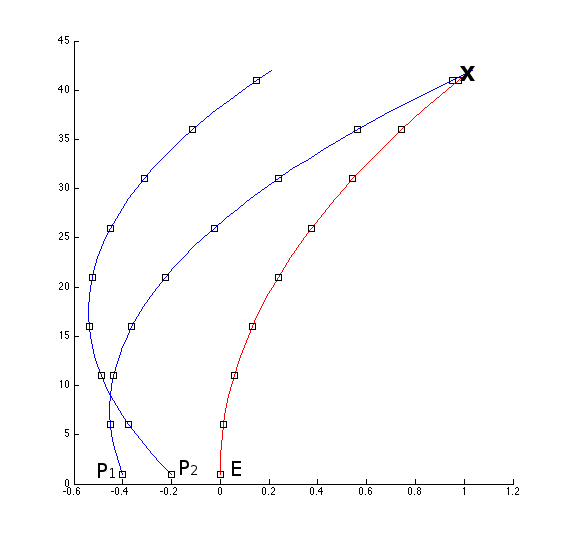}
 \caption{Optimal trajectories of the agents in the first component over time. $X$ denote the point of capture} \label{fig11}
\end{center}
\end{figure}
\noindent
For the special case when $m=2$, $d(x)=x$, $\alpha = 1$ and $\beta_{1}=\beta_{2}= 0.5$,  (Figure \ref{fig10}) shows computations of the value function with respect to the reduced coordinates 
$(y^1,y^2)=(x^1_1-x^2_1,x_2^1-x^2_2)$ in $\R^{2}$.
Figure \ref{fig11} shows an example of the evolution of the positions of the players over time,  with respect to the original coordinates.
Capture occurs at the point marked $X$, when pursuer $P_1$
overtakes the evader, despite starting farther from the evader than pursuer $P_{2}$.

%
\subsection{A Pursuit/Evasion Game With No Decomposition} 
\label{section4_3}
We now provide a simple example illustrating  that, for a multiple pursuers/single evader game, with target a union of target subsets, each associated with the evader and just one of the pursuers, may {\it fail} to  have a decomposition. In this example, it is possible to derive formulae for the value functions involved, and to test the conditions for decomposition directly. 
\ \\

\noindent
We denote by $(P^{3})$ the special case of $(P)$ in which $m_{1}=1$, $m_{2}=2$ and $n=1$.
$$
f_{1}(x_{1},a)=a,\; f_{2}(x_{2},b_{1})=b_{1} \mbox{ and }  f_{3}(x_{3},b_{2})=b_{2}\;.
$$
The controls actions of the players are constrained as follows:
$$
a \in A:= [-\alpha,+\alpha], 
$$
$$ b_{1} \in B_{1}:=[-1,+1]\,\mbox{ and }\, b_{2} \in B_{2}:=[-1,+1],
$$
for some $\alpha \in (0, 1)$.
We take the target to be 
$$
{\mathcal{T}} \,=\, {\mathcal{T}}_{2} \cup {\mathcal{T}}_{3}, \quad where
$$
$$
{\mathcal{T}}_{2}= \{(x_{1}, x_{2},x_{3})\,|\, x_{1}=x_{2} \}, \,
{\mathcal{T}}_{3}= \{(x_{1}, x_{2},x_{3})\,|\, x_{1}=x_{3} \}\, .
$$
(In this version of the game, two pursuers chase a single evader in $1D$ space. The game terminates when either pursuer meets the evader.)   Denote by $(P^{3}_{2})$ and $(P^{3}_{3})$ the modified games in which the target ${\mathcal{T}}$ is replaced by the subsets ${\mathcal{T}}_{2}$ and ${\mathcal{T}}_{3}$ respectively.
The HJI equation is
\begin{equation}
\label{HJE3}
 F^{3}(D_{x_{1}}u, D_{x_{2}}u, D_{x_{3}}u  )\,=\,0\, ,
\end{equation}
in which 
$$
F^{3}(p_{1},p_{2},p_{3} )\,=   |p_{2}|+ |p_{3}|- \alpha |p_{1}| -1\; .
$$

\noindent
Optimal strategies for  both games $(P^{3}_{2})$ and $(P^{3}_{3})$ are: the evader moves away from the pursuer, and the  pursuer moves towards the evader, as quickly as possible. A simple calculation based on these observations yields upper values for $(P^{3}_{2})$ and $(P^{3}_{3})$, namely:
\begin{eqnarray*}
&& u_{2}(x_{1}, x_{2}, x_{3})=(1-\alpha)^{-1}|x_{2}-x_{1}| \, ,
\\
&& u_{3}(x_{1}, x_{2}, x_{3})=(1-\alpha)^{-1}|x_{3}-x_{1}| \, ,
\end{eqnarray*}

for all $x=(x_{1}, x_{2}, x_{3})\in \R^{3}$.
Define $\bar u(.): \R^{3}\rightarrow \R$ to be
\begin{equation*}
\bar u(x) = \mbox{min} \{u_{1}(x),u_{2}(x)\}\,\mbox{ for } x \in \R^{3}\;.
\end{equation*}
\begin{proposition}
\label{section3_3}
$\bar u(.)$ is not a continuous viscosity solution for (\ref{HJE3}) on $\R^{3} \backslash {\mathcal{T}}$. 
\end{proposition}
Since the upper value for $(P^{3})$ is a viscosity solution on $\R^{3}\backslash {\mathcal T}$, vanishing on  ${\mathcal T}$, we may conclude that $\bar u(.)$ is not the value function for $(P^{3})$.\par
\noindent
{\bf Proof.}  Take any $z>0$ and let 
$\bar x= (0,z,-z)$. Then  $\bar x \in \R^{3} \backslash ({\mathcal{T}}_{2} \cup {\mathcal{T}}_{3} )$. Also, $u_{2}(\bar x)=u_{3}(\bar x)$, and $u_{2}(.)$ and $u_{3}(.)$ are continuously differentiable at $\bar x$. From the formulae for the value functions we have
\begin{eqnarray*}
\nabla u_{2}(\bar z)&=& (-(1-\alpha)^{-1},(1-\alpha)^{-1},0), \;\mbox{ and } \\
\nabla u_{3}(\bar z)&=&  ((1-\alpha)^{-1},0,-(1-\alpha)^{-1})\,.
\end{eqnarray*}
Then, for any $\lambda  \in (0,1)$,
\begin{eqnarray*}
F^{3}(\lambda \nabla_{x}u_{2}(\bar x)+ (1-\lambda) \nabla_{x}u_{3} (\bar x))
&=&
\frac{\lambda}{(1-\beta)} +  \frac{1- \lambda}{(1-\beta)}
\\
-
\frac{\beta}{(1-\beta)} (-\lambda  + (1-\lambda))-1&=&\frac{2\lambda \beta}{(1-\beta)} \;>\; 0\,. 
\end{eqnarray*}
So condition $(E)$ is violated. Then, $\bar u(.)$ cannot be a continous viscosity solution, by Prop \ref{prop2.2}, part (b)(iii).  \par
\noindent
The true value function $u(.)$ for  $(P^{3})$ is expressed in terms of the subset: \begin{eqnarray*}&&
{\mathcal{D}}\,=\, \{(x_{1},x_{2},x_{3})\in \R^{3}\,|\, \mbox{sgn}\{x_{2}-x_{1}\}= - \mbox{sgn}\{x_{3}-x_{1}\} 
\\
&&\hspace{0.3 in} \mbox{ and }  \frac{1-\alpha}{1+ \alpha}|x_{3}-x_{1}| < |x_{2}-x_{1}|  <     \frac{1+\alpha}{1- \alpha}|x_{3}-x_{1}|  \}\, .
\end{eqnarray*}
It is
\begin{eqnarray*}
&&
 u(x_{1}, x_{2}, x_{3})\;=\;
\left\{
\begin{array}{l}
\frac{1}{1-\alpha} \min\{ |x_{2}-x_{1}|,|x_{3}-x_{1}|\}\\
\hspace{0.6 in}\mbox{for } (x_{1}, x_{2}, x_{3}) \in \R^{3} \backslash {\mathcal{D}}
\\
\frac{1}{2} (|x_{2}- x_{1}|+ |x_{3}-x_{1}|) \\
\hspace{0.6 in}\mbox{for } (x_{1}, x_{2}, x_{3}) \in  {\mathcal{D}}\, .
\end{array}
\right.
\end{eqnarray*}
We see that $u(.)$ coincides with $\mbox{min}\{ u_{1}(x), u_{2}(x)\}$, for $x \in \R^{3}\backslash {\mathcal{D}}$. But
$$
u(x) < \mbox{$\mbox{min}\{ u_{1}(x), u_{2}(x)\}$, for $x \in {\mathcal{D}}$}\, .
$$
(The value function is constructed according to the heuristic: each of the pursuers always travels at maximum speed towards the evader. if both pursuers are on the same side of the evader, the evader travels at maximum speed in the opposite direction until the evader is hit. If, on the other hand, the evader is between the two pursuers, the evader travels at maximum speed away from the closest pursuer until the two pursuers are equidistant. The evader then stops until the evader is reached. A check is then carried out that the value function is a continuous viscosity solution of (\ref{HJE}), has a continuous extension to ${\mathcal{T}}$ on which it vanishes, and which is therefore the upper value of the game.) 


\section*{Appendix: Proof of Lemma \ref{control}}
\begin{small}
Consider first the case $z=z'$. Fix $t > 0$. We examine the optimal control problem of
$$
\left\{\begin{array}{l}
\mbox{Minimize $- y_{1}(t) $ }
\mbox{subject to }\\
(\dot y_{1}(s),\dot y_{2}(s))= (y_{2}(s), -d(y_{2}(s))+a(s)),\,\mbox{a.e. } s\in [0,t] \,,
\\ a(s)\in [-1,+1], \; \mbox{a.e. } s\in [0,t], \;
(y_{1}(0),y_{2}(0))= (z_{1},z_{2})\,.
\end{array}
\right.
$$
(Notice that the controlled differential equation in this problem is that governing the motion of the evader.)  The data for the problem satisfy standard hypotheses for the existence of a minimizer $a^{*}(.)$ on $[0,t]$, with corresponding state trajectory $y^{*}(.)$  (see, e.g. \cite[Chap. 2]{V00}). We can establish, by means of a simple contradiction argument, that the nonsmooth Maximum Principle (see \cite[Thm. 6.2.3]{V00}) applies in normal form. We deduce the existence of a costate arc $p(.)=(p_{1}(.),p_{2}(.))$ such that $p_{1}(.)\equiv +1$, and $p_{2}(.)$ satisfies the differential equation and right endpoint boundary condition
$$
-\dot p_{2}(s) =  +p_{1}(s)- \xi(s)\,p_{2}(s)\; \mbox{for } s \in [0,t] \mbox{  and } p_{2}(t)=0\,.
$$
 Here, $\xi_{i}(.)$ is a Lipschitz continuous function satisfying
$
\xi(s) \in \mbox{co}\, \partial_{L} \, d_{1}(y^{*}_{2}(s)) \,\mbox{a.e.}\,,
 $ in which $\partial_{L} \, d_{1}$ is the limiting subdifferential. The solution $p_{2}(.)$  is strictly positive on $[0,t)$. From the `maximization of the  Hamiltonian'  
$$
a^{*}(s) = \mbox{arg max}\, \{ p_{2}(s) a \,|\, a \in [-1,+1]  \} = +1 \, ,
$$  
$a^{*}(.)= \bar a(.)$ on $[0,t]$. 
We have shown that,  for any $t\geq 0$ and initial condition $z$, $a(.)=\bar a(.)$ maximizes $y_{1}(t)$. This confirms (\ref{monotone}) when $z'=z$.\\
\noindent
 We now show that (\ref{monotone}) is true also when $z^{'}_{1}=z_{1} $ and $z^{'}_{2} > z^{'}_{2}$. In view of the preceding analysis, we can assume that $a(.)= \bar a(.)$. Write $(y_{1}(.),y_{2}(.))$ and $(y^{'}_{1}(.),y^{'}_{2}(.))$ for the solutions to the state equation, for initial states $z= (z_{1},z_{2})$ and $z'=(z^{'}_{1}, z^{'}_{2})$ respectively. Take any time $\bar t>0$. By assumption $\dot y^{'}_{2}(0) > \dot y_{2}(0)$. So there are two cases to consider\\
(a): $\dot y^{'}_{2}(t) > \dot y_{2}(t)$ for all $t \geq 0$. In this case, since $y^{'}_{1}(0)-y_{1}(0) >0$, we have, as required,
$$
y^{'}_{1}(\bar t)- y_{1}(\bar t) = (y^{'}_{1}(0)-y_{1}(0) )+ \int _{0}^{\bar t} 
(\dot y^{'}_{2}(t) - \dot y_{2}(t))dt > 0\,.
$$
(b) There exists $t' \in (0, \bar t]$ such that $\dot y^{'}_{2}(t) > \dot y_{2}(t)$ for $t \in [0, t')$ and  $\dot y^{'}_{2}(t')= \dot y_{2}(t')$. In this case we show as, in the previous case, that $y^{'}_{1}(t')- y_{1}(t')>0$. We deduce from the  uniquess of solutions to the differential equation
$$
\dot y_{2}(t)= d(y_{2})+1\, ,
$$
on $[t', \bar t]$, for fixed initial condition, that $y^{'}_{2}(t) = y_{2}(t)$ for $t \in [t',\bar t]$. Hence, again, the required relation
$$
y^{'}_{1}(\bar t)- y_{1}(\bar t) = (y^{'}_{1}(t')-y_{1}(t') )+ \int _{t'}^{\bar t} 
(\dot y^{'}_{2}(t) - \dot y_{2}(t))dt > 0\,.
$$
\end{small}

\section*{Acknowledgements}
This  work has been supported by the European Union under the 7th Framework Programme FP7-PEOPLE-2010-ITN SADCO, "Sensitivity Analysis for Deterministic Controller Design".


\begin{thebibliography}{00}
\bibitem{BCD97} M. Bardi and I. Capuzzo-Dolcetta, \emph{Optimal Control and Viscosity Solutions of Hamilton-Jacobi-Bellman Equations}, Birkh\"auser, Boston 1997.
\bibitem{BRP99}M. Bardi, T.E.S. Raghavan, T. Parthasarathy, \emph{Stochastic and Differential Games: Theory and Numerical Methods}, Birkh\"auser, Boston, 1999.
\bibitem{B98}  G. Barles, \emph{Solutions de viscosit\`e des equations d'Hamilton--Jacobi}, Springer--Verlag, 1998.
\bibitem{CS04}  P. Cannarsa and C. Sinestrari, \emph{Semiconcave functions, Hamilton-Jacobi equations, and optimal control}, Birkh\"auser, Boston, 2004.
\bibitem{C89}W. Chodun, \emph{Differential games of evasion with many pursuers}, J. Math. Anal. Appl., 142(2) (1989) pp. 370--389. 
\bibitem{EK72} R.J. Elliott and N.J. Kalton,  \emph{Values in differential games}, Bull. Amer. Math. Soc.,  78(3) (1972) pp. 427--431.
\bibitem{CDCpaper} A. Festa and R.B. Vinter, \emph{A decomposition technique for  pursuit evasion games with many pursuers}, Proceedings of {52nd IEEE Control and Decision Conference (CDC)}, (2013), pp. 5797--5802.
\bibitem{F71}A. Friedman,  \emph{Differential Games},  John Wiley \& Sons, New York, USA, 1971.
\bibitem{I05}G.I. Ibragimov, \emph{Optimal pursuit of an evader by countably many pursuers}, Differ. Equ., 41(5) (2005) pp. 627–-635.
\bibitem{I65}R. Isaacs,  \emph{Differential Games}, John Wiley \& Sons, New York, USA, 1965.
\bibitem{IL81}R.P. Ivanov and Yu. S. Ledyaev, \emph{Time optimality for the pursuit of several objects with simple motion in a differential game},  Trudy Mat. Inst. Steklov., 158 (1981) pp. 87--97. 
\bibitem{KS88}N.N. Krasovskii and A.I. Subbotin,  \emph{Game-Theoretical Control Problems}, Springer, New York, 1988.
\bibitem{LT12}Y. Ledyaev and J.S. Treiman, \emph{Sub-and supergradients of envelopes, semicontinuous closures, and limits of sequences of functions}, Russ. Math. Surv., 67(2) (2012) pp. 345--373.
\bibitem{P76} B.N. Pshenichnii, \emph{Simple pursuit by several objects},  Cybern.  Syst. Anal., 12(3)  (1976) pp. 484--485.
\bibitem{VSH02}R. Vidal, O. Shakernia, J. Kim, H. Shim and S. Sastry, \emph{Probabilistic Pursuit--Evasion Games: Theory, Implementation, and Experimental Evaluation},  IEEE T. Robotic. Autom., 18(5) (2002) pp. 662--669.
\bibitem{V00}R. Vinter, \emph{Optimal Control}, Birkh\"auser, Boston, 2000.
\end{thebibliography}
\end{document}